\documentclass[10pt]{article}
\usepackage{amsmath,amssymb,amsthm,latexsym}

\newtheorem{theorem}{Theorem}
\newtheorem{prop}{Proposition}

\numberwithin{equation}{section}

\begin{document}
\title{{\sc Wrapping Brownian motion and heat kernels on compact Lie groups}}
\author{David Maher}
\date{}
\maketitle

\newcommand{\g}{\mathfrak{g}}
\newcommand{\kg}{\mathfrak{k}}
\newcommand{\tg}{\mathfrak{t}}
\newcommand{\ug}{\mathfrak{u}}
\newcommand{\p}{\mathfrak{p}}
\newcommand{\X}{\mathfrak{X}}
\newcommand{\af}{\mathfrak{a}}
\newcommand{\R}{\mathbb{R}}
\newcommand{\C}{\mathbb{C}}
\newcommand{\Z}{\mathbb{Z}}
\newcommand{\N}{\mathbb{N}}
\newcommand{\E}{\mathbb{E}}
\newcommand{\bbP}{\mathbb{P}}
\newcommand{\T}{\mathbb{T}}
\newcommand{\F}{\mathcal{F}}
\newcommand{\A}{\mathcal{A}}
\newcommand{\x}{\mathbf{x}}
\newcommand{\y}{\mathbf{y}}
\newcommand{\ad}{\mathrm{ad}}
\newcommand{\Ad}{\mathrm{Ad}}
\newcommand{\Exp}{\mathrm{Exp}}
\newcommand{\grad}{\mathrm{grad} \,}



\section{Introduction}

The partial differential equation given on $\R^n$ by

\begin{equation}\label{311}
\partial_t u(x, t) = \tfrac{1}{2} \Delta u(x, t), \phantom{abcde} t \in \R^+, \; x \in \R^n,
\end{equation}
where $\Delta$ is the Laplacian, represents the dissipation of heat
over a certain time.  The fundamental solution of the associated
semigroup $e^{t \Delta / 2}$, known as the {\bf heat kernel}, $p_t$
is given by a unique, strongly continuous, contraction semigroup of
convolution operators which may be convolved with the initial data
$f(x) = u(0,x)$ to give the solution to the Cauchy problem.  That
is,
$$
u(x,t) = e^{t \Delta / 2} f(x) = (p_t * f)(x) = \int_{\R^n} p_t(x-y) f(y) dy
$$

The heat kernel may also be expressed as the transition density of a {\bf Brownian motion}, $B_t$:
$$
p_t (x) = \E(B_t), \phantom{abcde} \text{moreover,} \; \; (p_t * f)(x) = \E(f(B_t))
$$

Similar statements hold when $\R^n$ is replaced by a Lie group.\\

In this article we will briefly demonstrate how these results may be
transferred from the Lie algebra (regarded as $\R^n$) to a compact
Lie group using the so-called wrapping map (\cite{DW2}).
Additionally, we shall provide the mechanism that allows one to
``wrap" a Brownian motion, and then find the heat kernel by taking
the expectation of the ``wrapped" process and applying a
Feynman-Ka\v{c} type transform.  We will also briefly discuss how
these results may be extended to compact symmetric spaces and
complex Lie groups.  Full details and proofs can be found in
\cite{M}.

\section{The wrapping map}

The wrapping map was devised by Dooley and Wildberger in \cite{DW2}.
 Let $G$ be a compact semisimple Lie group with Lie algebra $\g$.  We
define the {\bf wrapping map}, $\Phi$ by
\begin{equation}
\langle \Phi(\nu), f \rangle = \langle \nu,j \tilde{f} \rangle
\end{equation}
where $f \in C^\infty (G)$, $\tilde{f} = f \circ \exp$ and $j$ the
analytic square root of the determinant of the exponential map.  We
need to place some conditions on $\nu$ for $\Phi(\nu)$ to be
well-defined  -  this is the case when $\nu$ is a distribution of
compact support on $\g$, or $j\nu \in L^1(\g)$ .  We call
$\Phi(\nu)$ the \textit{wrap} of $\nu$.  The principal result is the
{\bf wrapping formula}, given by
\begin{equation}\label{wrapform}
\Phi(\mu *_\g \nu) = \Phi(\mu) *_G \Phi(\nu)
\end{equation}

This formula originated from their previous work on sums of adjoint orbits (\cite{DRW}), and can be considered as a global version of the Duflo isomorphism (\cite{DUF}).  The proof of (\ref{wrapform}) is particularly elegant, using only the Kirillov character formula and some abelian Fourier analysis.  Full details are in \cite{DW2}.\\

What (\ref{wrapform}) shows us is that problems of convolution of central measures or distributions on a (non-abelian) compact Lie group can be transferred to Euclidean convolution of Ad-invariant distributions on $\g$.\\

Thus, since the solution to the Cauchy problem for the heat equation can be written as a convolution between the heat kernel and the initial data, we should be able to wrap the heat kernel on $\g \cong \R^n$ to that on $G$, and transfer the corresponding solution of the Cauchy problem.\\

Given the remarks in section 1, it is clearly of interest also to consider whether there is a way to wrap Brownian motion to obtain the heat kernel on $G$.

\section{The wrap of Brownian motion}

Critical to wrapping a Brownian motion and the heat kernel from $\g$ to $G$ is how the infinitesimial generator of the respective process and semigroup  -  the Laplacian  -  is affected by wrapping.  The Laplacian on $\g$ is not quite wrapped to the Laplacian on $G$  -  a quantity that may be interpreted as a ``curvature" term arises.  More precisely, we have:
\begin{prop}\label{wraplap} Let $G$ be a compact connected Lie group with Lie algebra $\g$.  Then for any Schwartz function, $\mu$ on $\g$
$$
\Phi \bigl( L_\g (\mu) \bigr) = (L_G + \| \rho \|^2) \bigl( \Phi \mu \bigr)
$$
where $\Phi$ is the wrapping map, $L_\g$ is the Laplacian on $\g$ (regarded as a Euclidean vector space), $\rho$ the half sum of positive roots, and $\| \cdot \|$ the norm given by the Killing form.\\
\end{prop}

$L_G + \| \rho \|^2$ is also known as the {\bf shifted Laplacian}.  We shall refer the process and semigroup generated by $L_G + \| \rho \|^2$ as a {\bf shifted Brownian motion} and a {\bf shifted heat kernel}, respectively.\\

The actual mechanics of \textit{wrapping Brownian motion} are not immediately obvious, since the natural objects for the wrapping map to act on are distributions.\\

The wrapping map is a homomorphism from the algebra of Ad-invariant distributions on $C^\infty (\g)$ to the algebra of central distributions on $C^\infty (G)$, defined by $\phi \mapsto \phi \iota$ where $\iota : f \mapsto j.f \circ \exp$.\\

We ``wrap Brownian motion'' in an analogous way by considering the mapping $\iota$ in the context of It\^o stochastic differential equations.\\

Very briefly, we may construct a Brownian motion $(\zeta_t)_{t \geq 0}$ on $\g$ (regarded as the Lie group $\R^n$) as the solution to the Stratonovich S.D.E.:
\begin{equation}
d\zeta_t = \sum_{i=1}^n \frac{\partial \zeta_t}{\partial x_i} \circ dB_t^{(i)}, \phantom{abcde} \zeta_0 = 0
\end{equation}

This is really just a shorthand for the ``full" It\^o S.D.E.:
\begin{equation}\label{Itoong}
h(\zeta_t) = h(0) + \sum_{i=1}^n \int_0^t \frac{\partial h}{\partial x_i} (\zeta_t) dB_t^{(i)} + \tfrac 12 \sum_{i=1}^n \int_0^t \frac{\partial^2 h}{\partial x_i^2} (\zeta_t) dt
\end{equation}
where $h \in C^\infty_0 (\R^n)$.  Likewise, we define our shifted
Brownian motion on $G$ as the solution to the S.D.E.:
\begin{equation}
d\xi_t = \sum_{i=1}^n X_i (\xi_t) \circ dB_t^{(i)} + \tfrac 12 \|\rho\|^2 \xi_t dt, \phantom{abcde} \xi_0 = e.
\end{equation}
where $\bigl( X_i \bigr)_{i=1}^n$ is an orthonomal basis of the Lie
algebra, or in ``full" form:
\begin{equation}\label{ItoonG}
f(\xi_t) = f(e) + \sum_{i=1}^n \int_0^t (X_i f)(\xi_t) dB_t^{(i)} + \tfrac 12 \sum_{i=1}^n \int_0^t (X_i^2 f)(\xi_t)dt + \tfrac 12 \|\rho\|^2 \int_0^t f(\xi_t) dt
\end{equation}
where $f \in C^\infty (G)$.  To ``wrap of Brownian motion" we
replace $f \in C^\infty (G)$ with $j.f \circ \exp \in C^\infty_c
(\g)$, and let  $j.f \circ \exp = h \in C^\infty_0 (\g)$.  This can be shown to be
$$
h(\zeta_t) = h(0) + \sum_{i=1}^n \int_0^t \frac{\partial h}{\partial x_i}(\zeta_s) dB_t^{(i)} + \tfrac 12 \sum_{i=1}^n \int_0^t \frac{\partial^2 h}{\partial x_i^2}(\zeta_s)ds
$$
which is (\ref{Itoong}).  Thus we have
\begin{prop}\label{wrapbm}  Let $\zeta_t$ be a Brownian motion on $\g \cong \R^n$.  The wrap of $\zeta_t$ is a Brownian motion on $G$ with a potential of $\| \rho \|^2$, which we will call $\xi_t$.  That is,
$$
\Phi (\zeta_t) = \xi_t
$$
\end{prop}

We may now take expectations of each side to find the law of Brownian motion  -  the heat kernel  -  on $G$:

\begin{theorem}\label{wrapBM2}  Suppose $\xi_t$ is the wrap of the Brownian motion on $\g$, $\zeta_t$.  Then the law of $\xi_t$ may be found by wrapping the law of Brownian motion on its Lie algebra.  That is,
$$
\E_X(j. f \circ \exp(\zeta_t)) = \E_{\exp X} (f(\xi_t))
$$
which in law is given by
$$
\Phi (p_t)(\exp H) = q_t^\rho (g)
$$
where $p_t(x)$ is the heat kernel on $\g = \R^n$, and $q_t^\rho (g)$ is the heat kernel corresponding to the shifted Laplacian on $G$
\end{theorem}

The Feynman-Ka\v{c} theorem can be used to deal with the potential term $\| \rho \|^2$ to obtain a standard Brownian motion and heat kernel on $G$.  We omit the details, which will be presented in \cite{M}.

\section{The wrap of the heat kernel}

Let $p_t(x)$ be the heat kernel on $\R^n$, given by
\begin{equation}
p_t(x) = (2\pi t)^{-n/2} e^{-\frac{\| \x \|^2}{2t}}, \phantom{abcde} t \in \R^+, \; \x \in \R^n.
\end{equation}
and $q_t (g)$ is the heat kernel on $G$, given by
\begin{equation}
q_t (g) = \sum_{\lambda \in \Lambda^+ } d_\lambda \chi_\lambda (g) e^{-2\pi (\| \lambda + \rho \|^2 - \| \rho \|^2) t/2} , \phantom{abcde} t \in \R^+, \; g \in G.
\end{equation}
We write the shifted heat kernel on $G$ as $q_t^\rho (g)$, which is given by
\begin{equation}
q_t^\rho (g) = \sum_{\lambda \in \Lambda^+} d_\lambda \chi_\lambda (g) e^{-2\pi \| \lambda + \rho \|^2 t/2} , \phantom{abc} t \in \R^+, \; g \in G.
\end{equation}

Firstly, let's compute $\Phi(\nu)$.  When $\nu$ is suitably nice, it
has been shown in \cite{DW2} we can compute $\Phi(\nu)$ as a sum
over closed geodesics.  Let $\tg$ be the Lie algebra of the maximal torus, $T$, and let $\Gamma$ be the integer lattice in $\tg$, where $\Gamma = \{ H\in \tg \; : \; \exp (H) = e \}$.  We thus have:
\begin{equation}
\Phi(\nu) (\exp \, H) = \sum_{\gamma \in \Gamma} \Bigl(
\frac{\nu}{j} \Bigr) (H + \gamma), \; \; \; \forall H \in
\mathfrak{t}
\end{equation}
or secondly, as sum over highest weights $\Lambda^+$:
\begin{equation}
\Phi(\nu) (\exp \, H) = \sum_{\lambda \in \Lambda^+} d_\lambda
\nu^\wedge (\lambda + \rho) \chi_\lambda (g), \; \; \; \forall H \in
\mathfrak{t}
\end{equation}
which follows since it can be shown that $\Phi^\wedge (\nu) =
\nu^\wedge (\lambda + \rho)$ (see \cite{DW2}).  Equating these is
the Poisson summation formula for a compact Lie group.

From the above section, the law of $\xi_t$ may be found by wrapping the law of Brownian motion on its Lie algebra.  We put $p_t = \nu$ to find the law of the shifted Brownian motion on $G$:
\begin{align}
\Phi(p_t) (\exp H) & = \sum_{\lambda \in \Lambda^+} d_\lambda \, e^{-2\pi^2 \| \lambda + \rho \|^2 t} \chi_\lambda (H)\\
& = (2\pi t)^{-d/2} \sum_{n \in \Gamma} e^{\frac{-\|H + n\|^2}{2t}}
\frac{1}{j(H + n)}\label{4.7}
\end{align}
for all $H \in \tg$.  The first expression follows since $\hat{p}_t(\xi) = e^{-(2\pi \| \xi \|)^2t/2}$.

\section{Generalisations}

The wrapping formula needs some modification to hold for general
(compact) symmetric spaces $X$, equipped with tangent space $\p$,
with maximal abelian subalgebra $\af$.  This modification is
\begin{equation}\label{5.1}
\Phi(\mu *_{\p,e} \nu) = \Phi(\mu) *_X \Phi(\nu)
\end{equation}
where the convolution product on $\p$ is ``twisted" by a certain function $e$, which originates in the work of Rouvi\`ere \cite{R}.  See also \cite{D1}, \cite{D2}.\\

It is well-known in the physics literature that the ``sum over classical paths" does not hold
for general compact symmetric spaces (\cite{CAM}, \cite{DOW}).  That
is, performing a similar summation to (\ref{4.7}) to find the heat
kernel:
$$
\sum_{\gamma \in \Gamma^+} \Bigl( \frac{p_t}{j} \Bigr) (H + \gamma),
\phantom{abcde} \forall H \in \mathfrak{a}
$$
does not yield the (shifted) heat kernel on $X$.  The underlying
reason can be easily seen from (\ref{5.1}) in that we have a twisted
convolution on $\p$, which interferes with wrapping the heat
convolution semigroup:
$$
q_{t+s} = q_t *_X q_s = \Phi (p_t) *_{X} \Phi (p_s) = \Phi (p_t
*_{\p,e} p_s)
$$
which is not equal to $\Phi (p_{t+s})$.  It does turn out that we
can recover from this situation as the $e$-function and the
$j$-function are somewhat related. Basically, we need to consider
the heat kernel with potentials like $j^{-1} L_\p j$ on $\p$.  Even
for the 2-sphere this turns out to be difficult  - the potential in
this case is
$$
\frac{1}{H^2} - \text{cosec}^2 (H)
$$

We have also been able to extend our methods on wrapping Brownian
motion and heat kernels to some spaces where we know the the
wrapping formula holds.  A nice example are the complex Lie groups.
Instead of having to deal with a maximal torus $\T^n$, as in the
case of a compact Lie group, the subgroup corresponding to the
Cartan subalgebra is $(\R^+)^n$, so instead of summing over a
lattice, we just ``bend" the heat kernel from $\g$ to $G$ by
dividing by $j$, that is,
$$
\Phi (p_t) (\exp H) = (2\pi t)^{-n/2} \frac{1}{j(H)} \exp (-|H|^2 /2t), \phantom{abcde} H \in \af
$$

We can also wrap other processes - the key is to find how its
infinitesimal generator (call it $\mathcal{L_{\g}}$) wraps, that is,
$$
\Phi \bigl( \mathcal{L}_\g (u) \bigr) = (\mathcal{L}_G + C) \bigl( \Phi u \bigr)
$$

\section{Further directions}

\begin{list}{$\bullet$}{}
\item I am currently proving the wrapping formula for other Lie groups.  Once it is then known how to wrap a function, the heat kernel should then be able to be computed.  However, this is by no means straightforward  -  in the case of $SL(2,\R)$, the elements are conjugate to a choice of two abelian subgroups, isomorphic to $\T$ and $\R^+$.  Do we ``wrap" or ``bend"?  Probably both in some suitable fashion.
\item Wrapping the solutions of other P.D.E.'s.  In particular, any phenomena associated to them.  For example with the wave equation, what does it mean to ``wrap" Huygens' principle?  I should mention that it was for (odd dimensional) compact Lie groups, complex Lie groups, and the symmetric spaces $G / K$, $G$ complex, that Helgason was able to show that Huygens' principle holds when the shifted Laplacian is used (\cite{HE3}).
\item We would also like to know the $L^p - L^q$ bounds for a wrapped function.  For example, for what $p$ and $q$ do we have $ \| \Phi (u) \|_p \leq \| u \|_q$ ?  These could then be applied to obtain $L^p$ bounds of solutions of P.D.E.'s on Lie groups.  Currently, this is  only known when $p=q=1$.
\item These bounds could also be used to examine other behaviour such as convergence of Fourier transforms  -  if we used the ball multiplier, then in the case of compact Lie groups, our formula for $\Phi$ corresponds to the polygonal regions of positive weights considered for the convergence of Fourier series on compact Lie groups.
\end{list}

{\sc School of Mathematics, UNSW, Kensington 2052 NSW, Australia.}

Email: {\tt dmaher@maths.unsw.edu.au}

\end{document}